\documentclass[aop]{imsart}

\RequirePackage[OT1]{fontenc}
\RequirePackage{amsthm,amsmath}
\RequirePackage[numbers]{natbib}
\RequirePackage[colorlinks,citecolor=blue,urlcolor=blue]{hyperref}
\usepackage{amsmath}
\usepackage{amsthm}
\usepackage{amsfonts} 
\usepackage{hyperref}
\usepackage{epsfig}

\startlocaldefs
\numberwithin{equation}{section}
\theoremstyle{plain}

\endlocaldefs
\numberwithin{equation}{section}

\allowdisplaybreaks

\newtheorem{theorem}{Theorem}

\newtheorem{remark}{Remark}

\newtheorem*{ack}{Acknowledgement}

\renewcommand{\epsilon}{\varepsilon}

\newcommand{\1}[1]{{\mathbf 1}{\{#1\}}}

\newcommand{\T}{\mathbb{T}} 
\newcommand{\PR}{\mathbb{P}}
\newcommand{\Pb}{\overline{\mathbb{P}}}

\newcommand{\ES}{\mathbb{E}}

\begin{document}

\begin{frontmatter}
\title{Phase transition for the Once-reinforced random walk on $\mathbb{Z}^{d}$-like trees}
\runtitle{Phase transition for the Once-reinforced random walk}

\begin{aug}
\author{\fnms{Daniel} \snm{Kious}\thanksref{m3}\ead[label=e1]{daniel.kious@nyu.edu}}
\address{Daniel Kious\\ \printead{e1}}
\and
\author{\fnms{Vladas} \snm{Sidoravicius}\thanksref{m2,m3}\ead[label=e2]{vs1138@nyu.edu}}
\address{Vladas Sidoravicius\\ \printead{e2}}
\affiliation{Courant Institute of Mathematical Sciences, NYU, New York\thanksmark{m2}, and  NYU-ECNU Institute of Mathematical Sciences at NYU Shanghai\thanksmark{m3}}
\runauthor{D.~Kious and V.~Sidoravicius}
\end{aug}

\begin{abstract}
In this short paper, we consider the Once-reinforced random walk with reinforcement parameter $a$ on trees with bounded degree which are transient for the simple random walk. On each of these trees, we prove that there exists an explicit critical parameter $a_0$ such that the Once-reinforced random walk is almost surely recurrent if $a>a_0$ and almost surely transient if $a<a_0$.
This provides the first examples of phase transition for the Once-reinforced random walk.
\end{abstract}

\begin{keyword}[class=MSC]
\kwd[Primary ]{60K35}
\end{keyword}

\begin{keyword}
\kwd{Once-reinforced random walk}
\kwd{recurrence}
\kwd{transience}
\kwd{phase transition}
\end{keyword}

\end{frontmatter}

\section{Introduction}
We will be interested in non-Markovian random walks for which the future of the walk depends on its past trajectory. This fits in the large family of self-interacting random walks. Usually, these walks are hard to study and even basic properties such as recurrence and transience can be very difficult to obtain.\\
Here, we will focus on edge-interaction, where each edge of the considered graph has a current weight (depending on the past trajectory) and the walker jumps through an edge with a probability proportional to its weight.
One important example of such walks is the linearly Edge-reinforced random walk (ERRW) which was first introduced by Coppersmith and Diaconis \cite{CD}. The ERRW corresponds to the case where the current weight of an edge is equal to the number of time it has been crossed so far plus some initial weight $w$. In the eighties, it was conjectured that the ERRW on $\mathbb{Z}^d$ is recurrent if $d\le 2$ and, if $d\ge 3$, that it is recurrent for small $w$ and transient for large $w$. This remained a long standing open problem until the last few years. Let us state a brief history of the results obtained so far. We refer the reader to \cite{Kozsurvey,Pemsurvey,PTsurvey} for surveys on various random walks with reinforcement.\\
A phase transition was first proved on the binary tree by Pemantle \cite{Pemtree} who described the ERRW on trees as a random walk in independent random environment. Later, Merkl and Rolles considered, in a series of papers, the ERRW on various particular graphs, for instance proving recurrence on a graph which is $\mathbb{Z}^2$ with each edge replaced by $130$ edges in series, see \cite{MR}. One of the most recent and most important result is the  recurrence for $w$ small on $\mathbb{Z}^d$ which was proved independently by Angel, Crawford and Kozma \cite{ACK} and by Sabot and Tarr\`es \cite{ST}, at about the same time but with two different techniques. The proof in \cite{ACK} is intuitive and uses a nice simple idea, whereas the proof in \cite{ST} makes an explicit link with a model in quantum field theory, which seems to have far-reaching consequences. Using arguments from the physics litterature \cite{DSZ}, Disertori, Sabot and Tarr\`es \cite{DST} proved the transience for large $w$ on $\mathbb{Z}^d$, $d\ge3$. Up to now, the last striking result on ERRW is the recurrence on $\mathbb{Z}^2$ for any reinforcement parameter, proved by Sabot and Zeng \cite{SZ}. Research on ERRW is still going on: it seems that many links are to be discovered between the ERRW and models from physics or other interesting probabilistic models, as initiated by Sabot, Tarr\`es and co-authors \cite{ST,ST2,DST,STZ,SZ}.

Due to the absence of results about ERRW for many years, Davis \cite{Davis} introduced the Once-reinforced random walk (ORRW) as an \emph{a priori} simplified version of the ERRW. This is a walk for which the  current weight of an edge is $1$ if it has never been crossed and $a>1$ otherwise. It turned out that the study of this walk on $\mathbb{Z}^d$ is not easy at all and its recurrence/transience is still an open problem. It has been conjectured by the second author that the ORRW is recurrent on $\mathbb{Z}^d$ for  $d\in\{1,2\}$ and undergoes a phase transition for $d\ge3$, being recurrent when the parameter $a$ is large and transient when it is small. So far, there are only results on trees and particular graphs, for which no phase transition occurs. Sellke \cite{Sellke} proved that the ORRW is almost surely recurrent on the ladder $\mathbb{Z}\times\{1,...,d\}$ for $a\in(1,(d-1)/(d-2))$, see also \cite{Vervoort}. In contrast with the ERRW,  Durrett, Kesten and Limic \cite{DKL} showed that the ORRW is transient on regular trees for any $a>1$, which was later generalized to any supercritical tree by Collevecchio \cite{Coll}.\\
Until now, there was no example of graph on which the ORRW exhibits a phase transition and, among the results available so far, there is no good indication that a phase transition occurs. Here, we provide examples of irregular trees with bounded degree on which the ORRW indeed undergoes a phase transition. These trees have an overall drift that is similar to that of $\mathbb{Z}^d$. They were already introduced as a way to investigate the behavior of the simple random walk on $\mathbb{Z}^d$, see the survey \cite{DoyleSnell} by Doyle and Snell or Section \ref{sect_DS}.\\
While we do not claim that it really implies anything on $\mathbb{Z}^d$, it shows in a simple case that the once-reinforcement procedure can change the nature of the walk on graphs with a polynomial overall drifts.

\section{Construction of the trees and transience of simple random walk}

\begin{figure}[h]   
\includegraphics{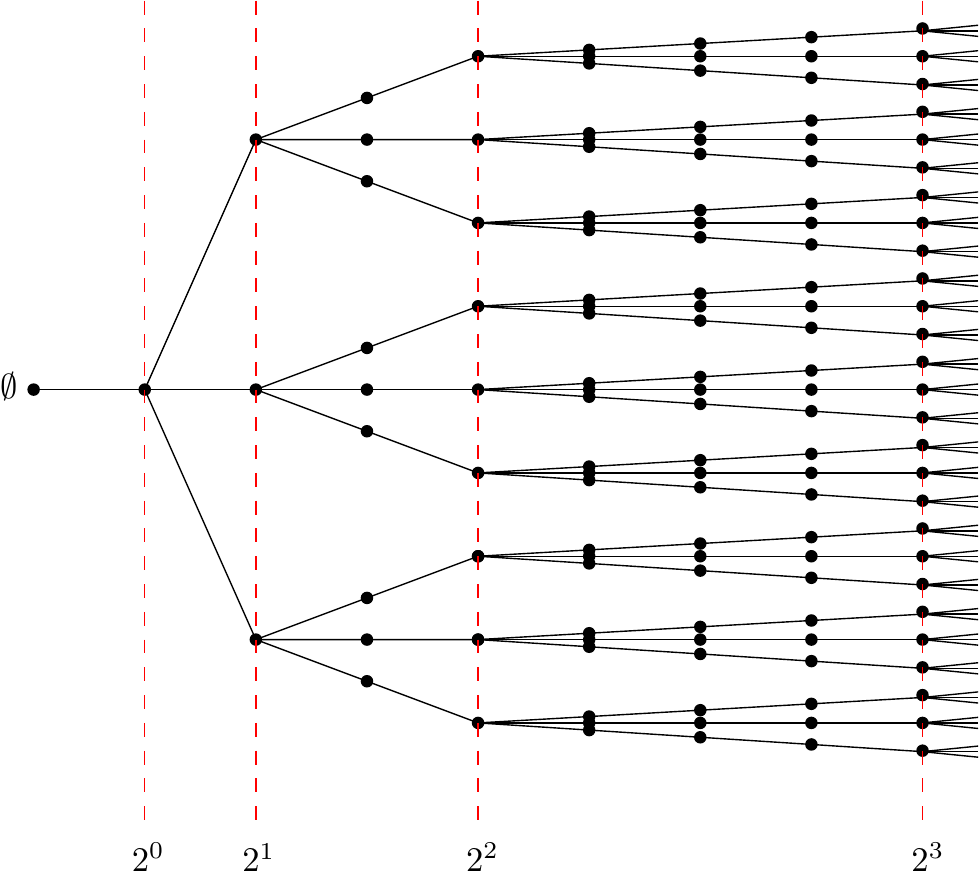}
   \caption{\label{arbre} The first generations of the tree $\T_3$.}
\end{figure}

Fix a positive integer $d\in\mathbb{N}\setminus\{0\}$ and let us construct a tree $\T_d$ rooted at a vertex $\emptyset$. The root is said to be at level $0$ and has one child (at level $1$). Each vertex at level $2^k$, $k\ge0$, has $d$ children and the vertices at any other level have only one child. This ``rarely splitting'' tree is depicted in Figure \ref{arbre}. For two vertices $x$, $y\in\T_d$, we write $x\sim y$ if they are neighbours, i.e.~if one them is the child of the other one, and we denote $[x,y]$ the non-oriented edge linking them. Besides, we call $x$ the ancestor of $y$ (at level $k$) if it is on the unique path between the root and $y$ (and if $x$ is at level $k$), and in this case $y$ is called a descendant of $x$.\\
Each edge of the tree has a conductance, or weight, $1$. For the simple random walk, this tree is known to be recurrent if $d\le 2$ and transient if $d\ge3$. Indeed, it is straightforward to compute the effective resistance from the root to infinity and obtain
\[
\mathcal{R}\left(\emptyset\leftrightarrow\infty\right)=1+\sum_{k=0}^{+\infty} \frac{2^k}{d^{k+1}},
\]
which is infinite for $d\in\{1,2\}$ and finite for $d\ge3$, see \cite{DoyleSnell} or \cite{LP} for more details on how to compute effective resistances and consequences on the recurrence/transience of the simple random walk.

\section{The Once-reinforced random walk on $\T_d$}
\subsection{Main result}

Let us first define the ORRW on $\T_d$, $d\ge1$, with parameter $a>0$.\\
The current weight of an edge is defined as follows: at time $n$, an edge has conductance $1$ if it has never been crossed (regardless of any orientation of the edges) and conductance $a>0$ otherwise. The case where $a>1$ matches the usual definition of the ORRW, but here we will also consider $a\le1$ which corresponds to some ``negative reinforcement'', or repulsion. For any $n\ge0$, let $E_n$ be the set of non-oriented edges crossed up to time $n$, that is
\[
E_n:=\left\{[x,y]: \ x,y\in\T_d\text{ and } \exists 1\le k\le n \text{ s.t.~}\{X_{k-1},X_k\}=\{x,y\}\right\}.
\]
At time $n\in\mathbb{N}$, if $X_n=x\in\T_d$, then the walk jumps to a neighbor $y\sim x$ with conditional probability
\[
\PR\left[\left.X_{n+1}=y\right|\mathcal{F}_n\right]=\frac{(a-1)\1{[x,y]\in E_n}+1}{\sum_{z:z\sim x}\left((a-1)\1{[x,z]\in E_n}+1\right)},
\]
where $\left(\mathcal{F}_n\right)$ is the natural filtration generated by the history of $(X_n)$.
\begin{theorem}\label{th}
For any positive integer $d$, the Once-reinforced random walk on $\T_d$ with reinforcement parameter $a>0$ is almost surely:
\begin{itemize}
\item[(i)] recurrent if $a>\log_2(d)$;
\item[(ii)] transient if $0<a<\log_2(d)$ .
\end{itemize}
\end{theorem}
\begin{remark}
We believe that the walk is recurrent when $a$ is critical for $d\ge2$, as it seems to corresponds to the criticality of some branching process, but we do not provide any proof. Nevertheless, for $d=2$, the critical case corresponds to the simple random walk on $\T_2$, which is recurrent. Also, note that the ORRW is always recurrent on the integer half-line $\T_1$ (obviously transient for the critical case $a=0$ which is degenerate) and on $\T_2$ if we restrict ourselves to the usual attractive case $a\ge1$.
\end{remark}

\subsection{The idea behind $\T_d$} \label{sect_DS}

The trees $\T_d$ were already introduced by Doyle and Snell \cite{DoyleSnell} in order to investigate the behavior of the simple random walk on $\mathbb{Z}^d$. The idea is pretty simple: we want a tree with spheres having a polynomial number of points with respect to their radius. Namely, in order to compare it with $\mathbb{Z}^d$, we want that, when we double the radius of a sphere, this multiplies the numbers of vertices on it by $2^{d-1}$. With this in mind, it is quite easy to see that $\T_2$ ``corresponds'' to $\mathbb{Z}^2$, $\T_4$ to $\mathbb{Z}^3$ and more generally $\T_{2^{d-1}}$ to $\mathbb{Z}^d$. In other words, for $d\ge1$, we think of $\T_d$ as a tree of dimension $\log_2(d)+1$.

In \cite{DoyleSnell}, the authors compute the effective resistance of these trees in order to embed them in $\mathbb{Z}^d$ and conclude about the recurrence or transience of the simple random walk. Although the phase transition of the ORRW on these trees is somehow good news for the conjecture of a phase transition, we do not claim that any strategy of embedding them on $\mathbb{Z}^d$ could work for the ORRW.


\subsection{Other examples of trees with phase transition}
We believe that the techniques used here can apply to various examples of trees with polynomial overall drift, in particular weighted regular trees.\\
Here is an example proposed to us by Gady Kozma. Consider the regular $3$-ary tree where the edges at level $k$ have a conductance $2^{-k}$. This tree is transient for the simple random walk (on each site, the local drift is positive). Applying the very same techniques as in this paper, one can prove that the ORRW on this tree is recurrent if $a>2$ and transient if $0<a<2$.

\subsection{The behavior of the ERRW on $\T_d$}
It is interesting to note the different behaviors of the ERRW and the ORRW on the trees $\T_d$. Indeed, while our main result reveals a phase transition for the ORRW, it can be proved that the ERRW is recurrent on $\T_d$, for any $d\ge1$ and for any reinforcement parameter (at least provided that all edges have the same initial weight).\\
We do not provide a full proof of this, as this does not provide anything new and it is fairly easy to obtain this using the results of \cite{Pemtree,LyoPem}. Let us give a rough blueprint. First, one can prove that, for any $\epsilon>0$, $\T_d$ is, with positive probability, a subtree of a Galton-Watson tree with mean offspring $1+\epsilon$, using for instance Lemma $6$ of \cite{Pemtree}. Second, one can prove that, for any reinforcement parameter, there exists $\epsilon>0$ such that the ERRW is almost surely recurrent on this Galton-Watson tree, using Theorem $3$ of \cite{LyoPem} and doing the same computation as in display $(6.2)$ of \cite{Pemtree} (note that the quantity computed there is strictly less than $1$ and does not depend on the precise law of the tree or on $\epsilon$).

\subsection{Strategy of the proof of Theorem \ref{th}}\label{sectstrat}
The strategy to prove recurrence is quite simple. First, we consider the ORRW on the half-line and estimate the probability for the walk to hit level $2^{k+1}$ before hitting the root once it is at level $2^k$. When $a$ is large enough, this probability is small compared to the number of points in the sphere of radius $2^{k+1}$ and this enables us to conclude recurrence.\\
The proof for transience is more subtle, adapting a very nice technique due to Collevecchio \cite{Coll}. The idea is to consider the walk after it has reached level $2^{kn_0}$, for some constant integer $n_0$, and observe the number $Z_k$ of children at level $2^{(k+1)n_0}$ that are hit before the walk goes back to the ancestor at level $2^{kn_0-1}$ (assuming that the walk eventually comes back to this ancestor). Then, we see the law of $Z_k$ as the offspring distribution at level $k$ of an inhomogeneous branching process. If $a$ is small enough, this branching process is supercritical and survives forever with lower-bounded probability. For the walk, it means that, for any $k_0$ and with lower-bounded probability, there exists an infinite path starting at level $2^{k_0n_0}$ such that the walk has to visit all the vertices on this path before going back to level $2^{k_0n_0-1}$. This easily implies the almost sure transience of the walk.

\subsection{Proof of Theorem \ref{th}}

\begin{proof}[Proof of $(i)$]
Fix $d\ge1$ and $a>\log_2(d)$, that is $a\ln(2)-\ln(d)>0$.\\Let us denote $L_n$, $n\ge0$, the set of vertices at level $n$ in $\T_d$. For $x\in\T_d$, define the following hitting time:
\[
T_x:=\inf\{n\ge0:X_n=x\},
\]
which is infinite if $x$ is never hit. Besides, for any $k\ge0$, define $T(2^k):=\inf_{x\in L_{2^k}}T_{x}$, the hitting time of level $2^k$. Also, for any $k\ge0$, define the first time the walker goes back the root after having hit level $2^k$, that is  $T_\emptyset^{(k)}=T(2^k)+T_\emptyset\circ\theta_{T(2^k)}$, where $\theta$ is the canonical shift. Note that $T(2^k)$ is almost surely finite for any $k\ge0$ as the walk cannot stay in any finite subtree forever. For any vertex $x\in L_{2^k}$, $k\ge1$, and for any $1\le n \le k$, we denote $\stackrel{\leftarrow}{x}_{{k-n}}\in L_{2^{k-n}}$ its unique ancestor at level $2^{k-n}$.\\
Notice that, for any $k\ge0$ and any integer $n\ge1$,
\begin{align}\nonumber
&\PR\left[\left.T_\emptyset^{(k)}>T(2^{k+n})\right| \mathcal{F}_{T(2^k)} \right]\\ \nonumber
&=\PR\left[\left.\bigcup_{x\in L_{2^{k+n}}}\left\{T_{\stackrel{\leftarrow}x_k}<\infty, T_{x}<T_\emptyset^{(k)},  T(2^{k+n})=T_x\right\}\right| \mathcal{F}_{T(2^k)}\right]\\ \label{circ1}
&\le \sum_{x\in L_{2^{k+n}}}\ES\left[\left.\1{T_{\stackrel{\leftarrow}x_k}<\infty}\PR\left[ \left. T_{x}\circ \theta_{T_{\stackrel{\leftarrow}x_k}}<T_\emptyset\circ \theta_{T_{\stackrel{\leftarrow}x_k}}\right|\mathcal{F}_{T_{\stackrel{\leftarrow}x_k}}\right]\right| \mathcal{F}_{T(2^k)}\right].
\end{align}

Now, we denote $X^1$, $T^1_\emptyset$ and $T^1(2^i)$ respectively the ORRW on the integer half-line (identifying the root $\emptyset$ and $0$), the hitting time of the root and the hitting time of level $2^i$ (or vertex $2^i$) associated to this walk. Also, we denote $\PR^1$ the probability measure associated to $X^1$. Then, for any $x\in L_{2^{k+n}}$, if $T_{\stackrel{\leftarrow}x_k}<\infty$, starting from $T_{\stackrel{\leftarrow}x_k}$ and up to $T_{x}\wedge T_\emptyset^{(k)}$, the jumps of $X$ along the unique path from the root to $x$ can be coupled with those of $X^1$ starting from $T^1(2^k)$. It is quite straightforward to define this coupling such that if $T_{x}\circ \theta_{T_{\stackrel{\leftarrow}x_k}}\wedge T_\emptyset\circ \theta_{T_{\stackrel{\leftarrow}x_k}}<\infty$ then $X_{T_{x}\circ \theta_{T_{\stackrel{\leftarrow}x_k}}\wedge T_\emptyset\circ \theta_{T_{\stackrel{\leftarrow}x_k}}}=X^1_{T^1(2^{k+n})\circ\theta_{ T^1(2^k) }\wedge T^1_\emptyset\circ\theta_{ T^1(2^k) }}$ almost surely, hence we will not detail it.\\
This yields
\begin{align}\nonumber
\1{T_{\stackrel{\leftarrow}x_k}<\infty}&\PR\left[ \left. T_{x}\circ \theta_{T_{\stackrel{\leftarrow}x_k}}<T_\emptyset\circ \theta_{T_{\stackrel{\leftarrow}x_k}}\right|\mathcal{F}_{T_{\stackrel{\leftarrow}x_k}}\right]\\ \label{circ2}
&\le\PR^1\left[  T^1(2^{k+n})<T^1(2^k)+T^1_\emptyset\circ \theta_{T^1(2^k)}\right].
\end{align}
In order to estimate this last quantity, note that when $X^1$ is at some level $j\in\{2^k,2^k+1,...,2^{k+n}-1\}$, with all the edges on its left having weight $a$ and all the edges on its right having weight $1$, the probability to hit $j+1$ before the root, at level $0$, is $j/(j+a)$ as can be seen by computing the effective resistance from $\emptyset$ to $j$ on the line or solving the Gambler's ruin problem, see \cite{LP}. Therefore, we have
\begin{align*}
&\PR^1\left[  T^1(2^{k+n})<T^1(2^k)+T^1_\emptyset\circ \theta_{T^1(2^k)}\right]\\
&\qquad\qquad=\prod_{j=2^k}^{2^{k+n}-1}\left(1-\frac{a}{j+a}\right)\le\exp\left(-a\sum_{j=2^k}^{2^{k+n}-1}\frac{1}{j+a}\right)\\
&\qquad\qquad\le \exp\left(-a\ln\left(\frac{2^{k+n}+a}{2^k+a}\right)\right)\le \exp\left(-an\ln(2)-a\ln\left(1-a2^{-k}\right)\right).
\end{align*}
Finally, notice that, in $\T_d$, there are $d^{k+n}$ vertices at level $2^{k+n}$, so that the last inequality, together with \eqref{circ1} and \eqref{circ2}, implies
\begin{align*}
&\PR\left[\left.T_\emptyset^{(k)}>T(2^{k+n})\right| \mathcal{F}_{T(2^k)}\right]\\
&\le\exp\left(k\ln(d)-n\left(a\ln(2)-\ln(d)\right)-a\ln\left(1-a2^{-k}\right)\right)
\end{align*}
and choosing $n=n_k=2k\lceil \ln(d+1)/\left(a\ln(2)-\ln(d)\right)\rceil$, which we can do because $a\ln(2)-\ln(d)>0$, we obtain, as $d+1\ge 2$, that
\begin{align*}
&\PR\left[\left.T_\emptyset^{(k)}>T\left(2^{k+n_k}\right)\right| \mathcal{F}_{T(2^k)}\right]\le\exp\left(-k\ln(2)-a\ln\left(1-a2^{-k}\right)\right).
\end{align*}
This last quantity is summable. Define the increasing sequence starting at $k_0=1$ and such that $k_i=k_{i-1}+2k_{i-1}\lceil \ln(d+1)/\left(a\ln(2)-\ln(d)\right)\rceil$, for any $i\ge1$. By Borel-Cantelli Lemma, there almost surely exists a finite random index $i_0$ such that, for any $i\ge i_0$, $T_\emptyset^{(k_i)}<T\left(2^{k_{i+1}}\right)$.
Hence, the root is almost surely visited infinitely many times, that is the walk is recurrent.
\end{proof}

\begin{proof}[Proof of $(ii)$]
Let us assume $a\ln(2)-\ln(d)<0$, with $d\ge2$ and $a>0$. We use a strategy which is a slight variation of the one in \cite{Coll}, as explained in Section \ref{sectstrat}. We will keep some notation from the proof of $(i)$.\\
Let ${n_0}\in\mathbb{N}$ the smallest integer such that, for any $k\ge1$,
\[
d^{n_0}\prod_{j=2^{k{n_0}}-2^{k{n_0}-1}}^{2^{(k+1){n_0}}-1}\left(1-\frac{a}{j+a}\right)\ge2.
\]
It can be easily checked that this integer ${n_0}$ exists and is finite as soon as $0<a<\ln(d)/\ln(2)$. It is important that ${n_0}$ does not depend on $k$.\\

For $k\ge 1$, fix a vertex $x \in L_{2^{k{n_0}}}$, denote $x_i$, $1\le i \le d^{n_0}$, its descendants at level $2^{(k+1){n_0}}$ and $x_{-1}$ its ancestor at level $2^{k{n_0}-1}$. We define $\T(x)$ the smallest connected subtree containing $x$, $x_{-1}$ and $x_i$, $1\le i \le d^{n_0}$, see Figure \ref{arbre3}.\\
\begin{figure}[h]   
\includegraphics{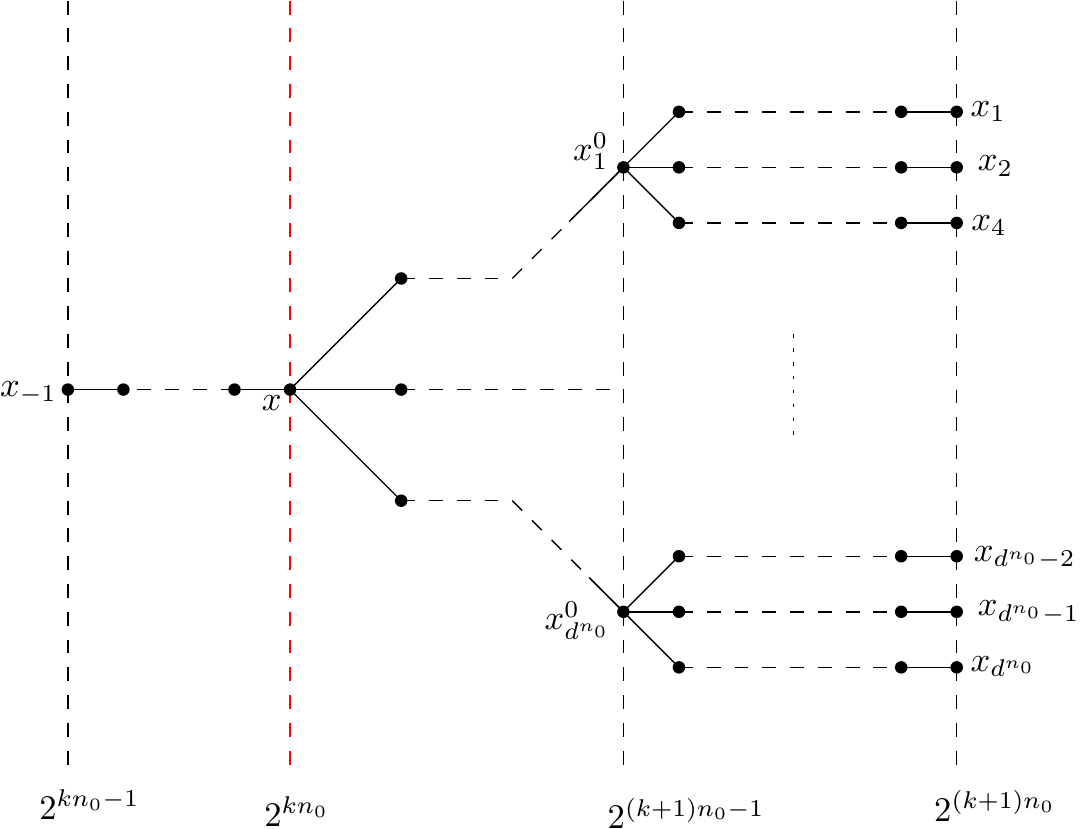}
   \caption{\label{arbre3} The finite tree $\T(x)$ and some distinguished points, in the case $d=3$.}
\end{figure}
Now, let us define $(Y^{x}_i)_i$ the ORRW on the finite tree $\T(x)$, reflected at $x_{-1}$ and $x_i$, $1\le i \le d^{n_0}$. We will denote $T^x_\cdot$ and $\PR^x$ respectively the stopping times and the probability measure associated to $Y^x$. We start $Y^{x}$ at $x$, putting weight $a$ on the edges on the path between $x_{-1}$ and $x$, and weight $1$ and the other edges. Then, consider this walk up to time $T^x_{ x_{-1}}$ and let $Z_k^x\in\{0,...,d^{n_0}\}$ be the number vertices among $x_i$, $1\le i \le d^{n_0}$, that have been visited before time $T^x_{x_{-1}}$. Note that the law of $Z_k^x$ depends on $k$, ${n_0}$, $a$ and $d$ but not on $x$ (as soon as it is at level $2^{kn_0}$).\\
For any $k\ge1$, $x\in L_{2^{kn_0}}$ and  $1\le i\le d^{n_0}$, by comparing $Y^x$ to a one-dimensional ORRW, we have that
\begin{align*}
&\PR^x\left[T^x_{x_i}<T^x_{x_{-1}}\right]=\prod_{j=2^{k{n_0}}-2^{k{n_0}-1}}^{2^{(k+1){n_0}}-1}\left(1-\frac{a}{j+a}\right),
\end{align*}
and thus
\[
\ES^x\left[Z_k^x\right]=\ES^x\left[\sum_{i=1}^{d^{n_0}}\1{T^x_{x_i}<T^x_{x_{-1}}}\right]=d^{n_0}\prod_{j=2^{k{n_0}}-2^{k{n_0}-1}}^{2^{(k+1){n_0}}-1}\left(1-\frac{a}{j+a}\right)\ge2.
\]
For $k_0\ge1$, let $B^{k_0}$ be the branching process where vertices at level $i$ have an offspring distribution given by that of $Z_{i+k_0}^x$ (for some $x\in L_{2^{(i+k_0)n_0}}$), and denote $B_i^{k_0}$ the number of particles at level $i$. Now, we want to use Theorem $1$ of \cite{BP} in order to prove that the process $B^{k_0}$ is supercritical. Let us precisely transcript the notation of \cite{BP} (using generating functions) according to our notation, for the reader's convenience. In Theorem 1 of \cite{BP}, $T$ denotes the extinction time of the branching process, thus 
$
P(T\le i)\equiv 1-\PR\left[B_i^{k_0}>0\right],
$
and we are therefore interested in the upper bound in display (2.4) of \cite{BP}. Besides, we have the following correspondence between the notation of \cite{BP} and our notation: $P_j\equiv \ES\left[B^{k_0}_j\right]$, $g_j'(1)\equiv\ES\left[Z_{k_0+j}^x\right]$ and $g_j''(1)\equiv{\ES\left[\left(Z_{k_0+j}^x\right)^2\right]-\ES\left[Z_{k_0+j}^x\right]}$, for any $x\in L_{2^{(j+k_0)n_0}}$ and for any $j\ge0$.\\
Hence, the upper bound in (2.4) of \cite{BP} and then using that
$\ES\left[\left(Z_k^x\right)^2\right]\le d^{2{n_0}}$,  $\ES\left[Z_k^x\right]\ge2$, for any $k\ge1$, and  $\ES\left[B^{k_0}_i\right]\ge2^i$, we have that
\begin{align*}
\lim_{i\to\infty}\PR\left[B_i^{k_0}>0\right]&\ge \limsup_{i\to\infty} \left[  \ES\left[B^{k_0}_i\right]^{-1}+\sum_{j=0}^{i-1}\frac{\ES\left[\left(Z_{k_0+j}^x\right)^2\right]-\ES\left[Z_{k_0+j}^x\right]}{\ES\left[Z_{k_0+j}^x\right]}    \ES\left[B^{k_0}_{j+1}\right]^{-1}  \right]^{-1}\\
&\ge \limsup_{i\to\infty}\frac{1}{2^{-i} + d^{2n_0}(1-2^{-i})} =d^{-2n_0}=:\delta
\end{align*}
where it should be noted that $\delta>0$ is a constant that does not depend on $k_0$.\\

Now, consider the following coloring scheme on the whole tree $\T_d$. Fix some $k_0\ge1$. For any $k\ge k_0$ and for any $x\in L_{2^{kn_0}}$, let the walks $Y^x$ defined above evolve independently up to time $T^x_{ x_{-1}}$. Pick any vertex $v_0\in L_{2^{k_0n_0}}$ and start by coloring it in white. By induction, for any $k\ge k_0$ and $x\in L_{2^{kn_0}}$, we color a descendant of $x$ at level $2^{(k+1)n_0}$ if $x$ is white and if this descendant has been visited by $Y^x$ before $x_{-1}$. Call the set of white points $B^{Y,k_0}$. It is straightforward to see that the set of white points is then equivalent to one realization of $B^{k_0}$.\\

The next step is to define a coupling of the trajectories $(Y^x)_{x\in\T_d}$ and $X$ such that, if $X$ comes back to level $2^{k_0n_0-1}$ after having hit level $2^{k_0n_0}$, then it has before visited all the points in  $B^{Y,k_0}$. Therefore, this can happen only if  $B^{Y,k_0}$ is finite and thus with probability at most $1-\delta$. The coupling is quite artificial but it hopefully clarifies that the dependencies created by the visits of $X$ in different generations of the tree are not too strong, or at least not relevant to prove a sufficient result. As the good coupling is not obvious, we detail it.\\

Fix two integers $k_0\ge1$, $K_0>k_0$ and wait for $X$ to hit level ${2^{k_0n_0}}$, denoting $v_0=X_{T({2^{k_0n_0}})}$ and  $v_{-1}$ the ancestor of $v_0$ at level ${2^{k_0n_0-1}}$. Let us define $V_{v_0}^{K_0}$ the set of vertices containing $v_0$ and all its descendants at levels ${2^{kn_0}}$, $k_0<k<K_0$. Besides, let us define $T_{\mathrm{end}}=T(2^{K_0n_0})\circ \theta_{T(2^{k_0n_0})}\wedge T({2^{k_0n_0-1}})\circ \theta_{T(2^{k_0n_0})}$. We are going to describe the evolution of a process
\[
W_n=\left( X_{T(2^{k_0n_0})+n\wedge T_{\mathrm{end}}}, C_n, \{(N_x(n),(Y^x_i)_{0\le i\le {N_x(n)}}), x\in V_{v_0}^{K_0}\}    \right),\]
where $C_n$ is a process taking values in $V_{v_0}^{K_0}\cup\{\emptyset\}$ starting at $C_0=v_0$, $Y^x_\cdot$ is a finite trajectory on the tree $\T(x)$ starting at $Y^x_0=x$ and $N_x(\cdot)$ is its associated clock (non-decreasing integer-valued process) starting at $N_x(0)=0$, for any $x\in V_{v_0}^{K_0}$. The process $C_n$ will indicate which of the $Y^x$ is currently coupled with $X$ (with the convention that $C_n=\emptyset$ means that there is currently no coupling). We denote $\Pb$ the probability measure associated to $(W_n)$. $(W_n)$ will be constructed in a way such that the trace (i.e.~the set of visited points) of $X$ on $\T(x)$ at time ${T(2^{k_0n_0})+n}$ matches the trace of $Y^x$ at time $N_x(n)$, for any $T_x\circ\theta_{T(2^{k_0n_0})}\le n \le T_{\mathrm{end}}$ (at these times, the edges on the left of $x$ are always considered visited for both walks).\\
We will need more notation.
For any $k\ge k_0$ and any $x\in L_{2^{k{n_0}}}$, we still denote $x_{-1}$ its ancestor at level $2^{kn_0-1}$ and $x_i$, $1\le i\le d^{n_0}$, its descendants at level $2^{(k+1)n_0}$, as well as $x_i^0$ the ancestor of $x_i$ at level $2^{(k+1)n_0-1}$ (they are not all distinct, see Figure \ref{arbre3}). We denote $H^x_n$ the (random) set of  vertices among $x_i$, $1\le i\le d^{n_0}$, visited by $X$ up to time ${T(2^{k_0n_0})+n\wedge T_{\mathrm{end}}}$ ($H^x_n$ can be empty). Also, define $\partial^x=\{x_{-1},x_1,...,x_{d^{n_0}}\}$, the boundary of the finite tree $\T(x)$.

%
Assume that we have constructed $W$ up to time $n<T_{\mathrm{end}}$ and that, for any $y\in V_{v_0}^{K_0}$ such that $T_y\le T(2^{k_0n_0})+n$, the traces of $X$ and $Y^y$ on $\T(y)$ are the same. Moreover, assume that $X_{{T(2^{k_0n_0})+n}}\in\T(x)\setminus\partial^x$ for some $x\in L_{2^{kn_0}}$, $k_0\le k<K_0$ and $C_n=x$. Then, let $X$ take one step according to its usual law, jumping on some vertex $z_{n+1}$ and:
\begin{itemize}
\item[$(i)$] If $n+1<T_{\mathrm{end}}$ then
\begin{itemize}
\item[$(a)$] If $z_{n+1}\notin \partial^x$, then $H^x_{n+1}=H^x_{n}$, $C_{n+1}=x$, $N_x(n+1)=N_x(n)+1$, $Y^x_{N_x(n+1)}=z_{n+1}$ and $N_v(n+1)=N_v(n)$ for any $v\in V_{v_0}^{K_0}\setminus \{x\}$.
\item[$(b)$] If $z_{n+1}=x_i\in \{x_1,...x_{d^{n_0}}\}\setminus H^x_{n}$ and thus $H^x_{n+1}=H^x_{n}\cup\{x_i\}$, then set $C_{n+1}=x_i$, $N_v(n+1)=N_v(n)$ for any $v\in V_{v_0}^{K_0}\setminus \{x\}$, set $Y^x_{N_x(n)+1}=x_i$ and let $Y^x$ evolve according to its usual law until it comes back to $x_i^0$, at some time $n'$, and set $N_x(n+1)=n'$. At the end of the construction, it will be clear that, in this case, $N_{x_i}(n+1)=0$ and that it is the first time $Y^{x_i}$ is coupled with $X$. The next time $X$ comes back to $x_i^0$, if it ever does, we will start again to couple it with $Y^x$ (see item $(d)$).
\item[$(c)$] If $z_{n+1}=x_i\in H^x_{n-1}$, then let $X$ move according to its usual law until it comes back to the father of $x_i$ (call this time $\widetilde{T}$ for the clock of $W$), setting $C_j=\emptyset$ for any $n+1\le j<\widetilde{T}\wedge T_{\mathrm{end}}$ and $C_{\tilde{T}\wedge T_{\mathrm{end}}}$ is $x$ if ${\widetilde{T}< T_{\mathrm{end}}}$ and $\emptyset$ otherwise. Besides, let $Y^x_{N_x(n)+1}=x_i$, $Y^x_{N_x(n)+2}=\stackrel{\leftarrow}{x}_i$ (the father of $x_i$) and, for any $n+1\le i\le \widetilde{T}\wedge T_{\mathrm{end}}$, set $N_x(i)=N_x(n)+2$ and $N_v(i)=N_v(n)$ for any $v\in V_{v_0}^{K_0}\setminus \{x\}$.
\item[$(d)$] If $z_{n+1}=x_{-1}$ then $N_x(n+1)=N_x(n)+1$, $Y^x_{N_x(n+1)}=x_{-1}$, $N_v(n+1)=N_v(n)$ for any $v\in V_{v_0}^{K_0}\setminus \{x\}$, and we set $C_{n+1}={\stackrel{\leftarrow}{x}_{{(k-1)n_0}}}$. From this time, the sub-tree starting at $x$ will not be relevant to us anymore, $Y^x$ has been constructed up to $T^x_{x_{-1}}$ and we will never couple it with $X$ again. The next time $X$ comes back to $x$, we will be in the situation $(c)$ described above.
\end{itemize}
\item[$(ii)$] If $n+1=T_{\mathrm{end}}$ then $C_{n+1}=\emptyset$, $N_x(n+1)=N_x(n)+1$, $Y^x_{N_x(n+1)}=z_{n+1}$ and $N_v(n+1)=N_v(n)$ for any $v\in V_{v_0}^{K_0}\setminus \{x\}$.
%
\end{itemize} 
The process $W$ is thus constructed up to time $T_{\mathrm{end}}$. Now, for any $x\in V_{v_0}^{K_0}$ such that $N_x(T_{\mathrm{end}})<T^x_{x_{-1}}$, we let $Y^x$ evolve according to its usual law up to time $T^x_{x_{-1}}=:N_x(T_{\mathrm{end}}+1)$, otherwise we already have $N_x(T_{\mathrm{end}})=T^x_{x_{-1}}=:N_x(T_{\mathrm{end}}+1)$. Finally, we set $N_x(i)=N_x(T_{\mathrm{end}}+1)$ for any $x\in V_{v_0}^{K_0}$ and $C_i=\emptyset$, for all  $i\ge T_{\mathrm{end}}+1$. This completes the construction of the process and we obtain
\[
W_\infty=\left( X_{T(2^{k_0n_0})+T_{\mathrm{end}}}, \emptyset, \{(T^x_{x_{-1}},(Y^x_i)_{0\le i\le T^x_{x_{-1}}}), x\in V_{v_0}^{K_0}\}    \right).
\]

The important point is to decouple $Y^{x}$ and $X$ when they hit for the first time some vertex $x_i$, $1\le i\le d^{n_0}$, and until $X$ comes back to $x_i^0$. This allows us to start coupling $X$ with $Y^{x_i}$ at time $T_{x_i}$. Besides, it is quite straightforward that the marginal laws of $X$ and $\{(Y_i^x)_{0\le i\le T^x_{x_{-1}}},x\in V_{v_0}^{K_0}\}$ under $\Pb$ are their usual laws (because the traces of $X$ and $Y^x$ match when they move together) and that the trajectories $(Y^x)_{x\in V_{v_0}^{K_0}}$ are independent (because the steps of $X$ when $C_n\neq x$ have no influence on the steps of $Y^x$).

Also, note that if $X_{T(2^{k_0n_0})+T_{\mathrm{end}}}=v_{-1}$, then we necessarily have that if $x\in V_{v_0}^{K_0}$ and $T_x<T(2^{k_0n_0})+T_{\mathrm{end}}$ then $N_x(T_{\mathrm{end}})=T^x_{x_{-1}}$ (i.e.~$X$ and $Y^x$ have been coupled until they hit $x_{-1}$ together). Thus, if we color the points in the tree $\bigcup_{x\in V_{v_0}^{K_0}} \T(x)$, starting with $v_0$ and defining $B^{Y,k_0}$ as above, it means that if $T_{v_{-1}}\circ \theta_{T_{v_0}}<T(2^{K_0n_0})\circ \theta_{T_{v_0}}$ then $X$ visits all the white points. Therefore, if $X_{T(2^{k_0n_0})+T_{\mathrm{end}}}=v_{-1}$ then $B^{Y,k_0}$ necessarily dies out before level $K_0-k_0$.
We thus obtain
%
\begin{align*}
&\PR\left[\left.T(2^{k_0n_0-1})\circ \theta_{T(2^{k_0n_0})}<\infty\right|\mathcal{F}_{T(n_0k_0)}  \right]\\
&=\lim_{K_0\to\infty} \PR\left[\left.T(2^{k_0n_0-1})\circ \theta_{T(2^{k_0n_0})}<T(2^{K_0n_0})\circ \theta_{T(2^{k_0n_0})}\right|\mathcal{F}_{T(n_0k_0)}  \right]\\
&\le\lim_{i\to\infty}\PR\left[B_i^{k_0}=0\right]\le 1-\delta.
\end{align*}
As $\delta>0$ does not depend on $k_0$, this implies that the walk is almost surely transient.
\end{proof}

\begin{ack}
The authors are very grateful to Pierre Tarr\`es for pointing out reference \cite{Coll}, which led to the proof for transience. Besides, they warmly thanks Gady Kozma and Bruno Schapira for motivating brain-teasing discussions on the topic. DK is grateful to the Ecole Polytechnique de Lausanne (EPFL) to which he was affiliated to at the time this work was completed.
\end{ack}


\begin{thebibliography}{99}
\bibitem{ACK} Angel, O., Crawford, N.~and Kozma, G.~(2014). Localization for linearly edge reinforced random walks. {\it Duke Math. J.} {\bf 163(5)}, 889--921.
\bibitem{BP}Agresti, A.~(1975). On the Extinction Times of Varying and Random Environment Branching Processes. {\it J. of Appl. Probab.} {\bf 12(1)}, 39--46.
 \bibitem{Coll} Collevecchio, A. (2006).  On the transience of processes defined on Galton-Watson trees. {\it Ann. Probab.} {\bf 34(3)}, 870--878.
  \bibitem{CD} Coppersmith, D.~and Diaconis, P.~(1986). { Random walks with reinforcement}. {\it Unpublished manuscript.}
 \bibitem{Davis} Davis, B. (1990). { Reinforced random walk}. {\it Probab.~Theory and Related Fields}, {\bf 84(2)}, 203--229.
 \bibitem{DST} Disertori, M., Sabot, C.~and Tarr\`es, P.~(2015). Transience of edge-reinforced random walk. {\it Comm. Math. Phys.} {\bf 339(1)}, 121--148.
 \bibitem{DSZ} Disertori, M., Spencer, T.~and Zirnbauer, M.~R.~(2010). Quasi-diffusion in a 3D supersymmetric hyperbolic sigma model. {\it Comm. Math. Phys.} {\bf 200(2)}, 435--486.
 \bibitem{DoyleSnell} Doyle, P.G. and Snell E.J. (1984). { Random walks and electrical networks}. {\it Carus Math. Monographs}, {\bf 22}, Math. Assoc. Amer., Washington, D.C.
 \bibitem{DKL} Durrett, R., Kesten, H.~and Limic, V.~(2002). Once edge-reinforced random walk on a tree. {\it Probab. Theory Related Fields}. {\bf 122(4)}, 567--592.
 \bibitem{Kozsurvey} Kozma, G.~(2012). Reinforced random walks. {\it Preprint}. arXiv:1208.0364.
  \bibitem{LyoPem} Lyons, R.~and Pemantle, R. (1992). Random Walk in a Random Environment and First-Passage Percolation on Trees. {\it Ann. Probab.} {\bf 20(1)}, 125--136.
 \bibitem{LP} Lyons, R. and Peres Y. (2016). Probability on trees and networks. Cambridge University Press. To appear,  \url{http://pages.iu.edu/\string~rdlyons/prbtree/prbtree.html}.
 \bibitem{MR} Merkl, F.~and Rolles, S.~W.~W.~(2009). Recurrence of the Edge-reinforced random walk on a two-dimensional graph. {\it Ann. Probab.} {\bf 37(5)}, 1679--1714.
  \bibitem{Pemtree} Pemantle, R. (1988). Phase transition in reinforced random walk and RWRE on trees. {\it Ann. Probab.} {\bf 16}, 1229--1241.
  \bibitem{Pemsurvey} Pemantle, R.~(2007). A survey of random processes with reinforcement. {\it Probab. Surv.} {\bf 4}, 1--79.
  \bibitem{ST} Sabot, C.~and Tarr\`es, P.~(2015). Edge-reinforced random walk, vertex-reinforced jump process and the supersymmetric hyperbolic sigma model. {\it J. Eur. Math. Soc.} {\bf 17(9)}, 2353--2378.
    \bibitem{ST2} Sabot, C.~and Tarr\`es, P.~(2015). Inverting Ray-Knight identity. {\it To appear in Probab.~Theory and Related Fields}.
        \bibitem{STZ} Sabot, C., Tarr\`es, P.~and Zeng, X.~(2015). The Vertex Reinforced Jump Process and a Random Schrödinger operator on finite graphs.{\it preprint}, 	arXiv:1507.04660.
  \bibitem{SZ} Sabot, C.~and Zeng, X.~(2015). A random Schrödinger operator associated with the Vertex Reinforced Jump Process on infinite graphs. {\it preprint}, 	arXiv:1507.07944.
  \bibitem{Sellke} Sellke, T.~(2006). Recurrence of reinforced random walk on a ladder. {\it Electron. J. Probab.} {\bf 11}, 301--310.
  \bibitem{PTsurvey} Tarr\`es, P.~(2011). Localization of reinforced random walks. {\it Preprint}. arXiv:1103.5536.
    \bibitem{Vervoort} Vervoort, M. (2002). {\it Reinforced random walks}. 
\end{thebibliography}
\end{document}